\documentclass{article}[15pt]
\usepackage{amsmath,amsmath, epsfig,amssymb, multicol}
\usepackage[active]{srcltx}
\usepackage{epstopdf}
\usepackage{graphicx}
\usepackage{color}
\usepackage{setspace}
\newtheorem{lemma}{Lemma}
\newtheorem{theorem}{Theorem}
\newtheorem{conjecture}{Conjecture}

\newtheorem{claim}{Claim}

\newcommand{\h}{{\cal  H}}
\newcommand{\ex}{{\rm ex}}
\newcommand{\C}{{\cal C}}

\begin{document}
\title{The number of edges in graphs with bounded clique number and circumference}

\author{{Chunyang Dou\footnote{School of Mathematical Sciences, Anhui University,
Hefei  230601, P.~R.~China. Email:{\tt chunyang@stu.ahu.edu.cn}. Supported in part by the National Natural Science Foundation of China
(No.\ 12071002, No.\ 12471319) and the Anhui Provincial Natural Science Foundation (No. 2208085J22).}}
~~~
{Bo Ning\footnote{ College of Computer Science, Nankai University, Tianjin, 300350, P.R. China. Emai: {\tt bo.ning@nankai.edu.cn}. Supported in part by the National Natural Science Foundation of China (No.\ 12371350) and the Fundamental Research Funds for the Central Universities (No.\ 63243151).}}
~~~
 {Xing Peng\footnote{Center for Pure Mathematics, School of Mathematical Sciences, Anhui University,
Hefei  230601, P.~R.~China. Email: {\tt x2peng@ahu.edu.cn}. Supported in part by the National Natural Science Foundation of China
(No.\ 12071002, No.\ 12471319) and the Anhui Provincial Natural Science Foundation (No. 2208085J22). } }}
\date{}
\maketitle

\begin{abstract}
   Let $\h$ be a family of graphs. The Tur\'an number $\ex(n,\h)$ is the maximum possible number of edges in an $n$-vertex graph which does not contain any  member  of $\h$ as a subgraph. As a common generalization of Tur\'an's theorem  and Erd\H{o}s-Gallai theorem on the Tur\'an number of matchings, Alon and Frankl determined    $\ex(n,\h)$ for $\h=\{K_r,M_k\}$, where $M_k$ is a matching of size $k$. Replacing $M_k$ by $P_k$, Katona and Xiao obtained the Tur\'an number of $\h=\{K_r,P_k\}$ for $r \leq \lfloor k/2 \rfloor$ and sufficiently large $n$. In addition, they proposed a conjecture for   the case where $r \geq \lfloor k/2 \rfloor+1$ and $n$ is  sufficiently large. Motivated by the fact that the result for $\ex(n,P_k)$ can be deduced from  the  one for  $\ex(n,\C_{\geq k})$, we investigate the Tur\'an number of $\{K_r, \C_{\geq k}\}$ in this paper. In other words, we aim to  determine the maximum number of edges in graphs with clique number at most $r-1$ and circumference at most $k-1$. For $\h=\{K_r, \C_{\geq k}\}$, we are able to show the value of $\ex(n,\h)$ for $r \geq \lfloor (k-1)/2\rfloor+2$ and all $n$. As an application of this result, we confirm Katona and Xiao's conjecture in a stronger form.  For $r \leq \lfloor (k-1)/2\rfloor+1$, we manage to show the value of $\ex(n,\h)$  for sufficiently large $n$.




\end{abstract}

\section{Introduction}\label{sec1}
The study of Tur\'an number of graphs is one of the central topics in extremal graph theory. For a family of graph $\h$, a graph is $\h$-free if it does not contain a member of $\h$ as a subgraph. The {\it Tur\'an number} $\ex(n,\h)$ is the maximum possible number of edges in an $n$-vertex $\h$-free graph. If $\h$ contains one graph $H$, then we  write  $\ex(n,H)$ instead of $\ex(n,\h)$. Mantel Theorem states that  $\ex(n,K_3)=\lfloor n^2/4\rfloor$ and the only extremal graph is the balanced bipartite graph. For cliques, the famous Tur\'an Theorem tells us  that $\ex(n,K_{p+1})=e(T(n,p))$ and the only extremal graph is $T(n,p)$, where $T(n,p)$ is the balanced complete $p$-partite graph with $n$ vertices. The number of edges in $T(n,p)$ is denoted by $t(n,p)$.
  Tur\'an's Theorem is viewed as the origin of extremal graph theory. For the matching $M_k$ of size $k$, Erd\H{o}s and Gallai \cite{EG59} determined the value of $\ex(n,M_k)$.
Let $P_k$ be the path with $k$ vertices and $\C_{\geq k}$ be the set of cycles of length at least $k$. In the same paper, Erd\H{o}s and Gallai proved the following results.
\begin{theorem}[Erd\H{o}s-Gallai \cite{EG59}]\label{EGpath}
$\ex(n,P_k) \leq \frac{(k-2)n}{2}$, where $k \geq 2$.
\end{theorem}

\begin{theorem}[Erd\H{o}s-Gallai \cite{EG59}]\label{EGcycle}
$\ex(n, \C_{\geq k}) \leq \frac{(k-1)(n-1)}{2}$, where $k \geq 3$.
\end{theorem}
If $k-1 | n$, then the union of  $n/(k-1)$ vertex disjoint $K_{k-1}$'s shows that   the upper bound in Theorem \ref{EGpath} is tight.
Similarly, if $k-2 |n-1$, then the graph consists of  $(n-1)/(k-2)$ copies of  $K_{k-1}$  sharing a common vertex gives the tightness of the upper bound in Theorem \ref{EGcycle}. Additionally, one can see Theorem \ref{EGpath} is a direct consequence of Theorem \ref{EGcycle}.

For general $n$, Faudree and Schelp \cite{FS75} and independently Kopylov \cite{K77} determined the  value of $\ex(n,P_k)$. Meanwhile, Woodall \cite{W76} obtained the value of $\ex(n,\C_{\geq k})$ and showed the result for $\ex(n,P_k)$ is a corollary of the one for $\ex(n,\C_{\geq k})$.

In the study of Tur\'an number of graphs, the celebrated Erd\H{o}s-Stone-Simonovits Theorem \cite{ES,ES1}  gives an asymptotic formula for $\ex(n,\h)$ when $\h$ contains no bipartite graphs:
\[
\ex(n, \h) = \left(1 - \frac{1}{p(\h)} \right) \binom{n}{2} + o(n^2),
\]
where $p(\h)=\min\{\chi(H)-1: H \in \h \}$.

Many existing results study the Tur\'an number of a single graph. For the case of $\h$ with two graphs, a classical result by Chv\'{a}tal and Hanson \cite{CH} gives the Tur\'an number of $\h=\{M_{\nu+1},K_{1,\Delta+1}\}$. A special case where $\nu=\Delta=k$ was first proved by Abbott, Hanson, and Sauer \cite{AHS}.
As a common generalization of Tur\'an's Theorem and Erd\H{o}s-Gallai Theorem for the Tur\'an number of matchings, Alon and Frankl \cite{AF} recently obtained the value of $\ex(n,\h)$ for $\h=\{K_r,M_k\}$. Wang, Hou, and Ma \cite{WHM} proved a spectral analogue of this result.   Taking a step forward, Katona and Xiao \cite{KX} investigated the Tur\'an number of  the family $\h=\{K_r,P_k\}$. To state  their results,  we need to define the join of two graphs.  Given two vertex disjoint graphs $F_1$ and $F_2$, the {\it join}  $F_1 \vee F_2$  is a graph obtained from $F_1 \cup F_2$ by connecting each vertex of $F_1$ and each vertex of $F_2$.

For $r \leq \lfloor k/2 \rfloor$, Katona and Xiao proved the following result.
\begin{theorem} [Katona-Xiao \cite{KX}]
Let $\h=\{K_r,P_k\}$. If $r \leq \lfloor k/2 \rfloor$ and $n$ is large enough, then
\[
\ex(n,\h)=\ex(\lfloor k/2 \rfloor -1,K_{r-1})+(\lfloor k/2 \rfloor -1)(n-\lfloor k/2 \rfloor +1),
\]
and $T(\lfloor k/2 \rfloor -1,r-2) \vee I_{n-\lfloor k/2 \rfloor +1}$ is an extremal graph.
\end{theorem}
For $r \geq \lfloor k/2 \rfloor+1$, they  made the following conjecture.
\begin{conjecture} [Katona-Xiao \cite{KX}] \label{conj1}
Assume $\lfloor k/2 \rfloor+1 \leq r < k$. If $k \geq 3$ is odd and $k-1 | n$, then $\ex(n,\h)=\tfrac{n}{k-1}\ex(k-1,K_r)$ for sufficiently large $n$. If $k \geq 4$ is even, then  $\ex(n,\h)=\binom{\lfloor k/2 \rfloor -1}{2}+(\lfloor k/2 \rfloor -1)(n-\lfloor k/2 \rfloor +1)$ for  sufficiently large $n$.
\end{conjecture}
 Katona and Xiao \cite{KX} also conjectured the asymptotic value of $\ex(n,\{P_k,H\})$ for a nonbipartite graph $H$. This conjecture was proved by Liu and Kang \cite{LK}.

As we mentioned above, the result for the Tur\'an number of $P_k$ is an easy corollary of the one for the Tur\'an number of $\C_{\geq k}$. In the same spirit,  we study the Tur\'an number of $\h=\{K_r, \C_{\geq k}\}$ in this paper. If $r \geq k$, then we can see $\ex(n,\h)=\ex(n, \C_{\geq k})$. Moreover, if $n \leq k-1$, then $\ex(n,\h)=\ex(n,K_r)$. Thus we assume
\begin{equation} \label{assump}
3 \leq r <k \leq n
\end{equation}
  throughout this paper.

  Assume that $n-1=p(k-2)+q$ with $q \leq k-3$. Define $F(n,k,r)$ as the graph which consists of $p$ copies of $T(k-1,r-1)$ and one copy of $T(q+1,r-1)$ sharing a common vertex. Let
\[
f(n,k,r)=p\ex(k-1,K_r)+\ex(q+1,K_r).
\]
Note that $f(n,k,r)=e(F(n,k,r))$.  Additionally, let
$$G_1=K_{\lfloor (k-1)/2 \rfloor} \vee I_{n-\lfloor (k-1)/2 \rfloor }$$
and
$$G_2=T(\lfloor (k-1)/2 \rfloor,r-2) \vee I_{n-\lfloor (k-1)/2 \rfloor }.$$
Our first result is the following one.
\begin{theorem} \label{thm1}
Let  $\h=\{K_r,\C_{\ge k}\}$ and $n-1=p(k-2)+q$. If $k\geq 6$ is even and  $\lfloor (k-1)/2 \rfloor+2\le r<k \leq n$, then
\[
\ex(n,\h)=f(n,k,r).
\]
 If $k \geq 5$ is odd and  $\lfloor (k-1)/2 \rfloor+2\le r<k \leq n$, then
\[
\ex(n,\h)=\max\left\{f(n,k,r), e(G_1)\right\}.
\]
\end{theorem}

\noindent
{\bf Remark 1:} For even $k$, if $k=4$,  then $r=3$ by the lower bound assumption for $r$. Furthermore, the assumption $G$ being $\{K_3,\C_{\geq 4}\}$-free implies that $G$ does not contain cycles. Thus $G$ either is a tree or is a forest. This case is trivial.   As  the assumption $3 \leq r<k$, it is natural to assume $k \geq 5$ for odd $k$.

 We prove the next theorem for $r\le\lfloor (k-1)/ 2\rfloor+1$.
\begin{theorem} \label{thm2}
Let  $\h=\{K_r,\C_{\ge k}\}$. If $r\le\lfloor (k-1)/ 2\rfloor+1$ and $n \ge k\ge 5$, then
\[
\ex(n,\h)=e(G_2)
\]
provided   $n \geq \tfrac{k^3}{4}$.
\end{theorem}
As an application of Theorem \ref{thm1}, we prove the following result.
\begin{theorem} \label{conj11}
Let $\h=\{K_r,P_k\}$ and $\lfloor k/2 \rfloor+1 \leq r <k \leq n$.  Assume that $n=p(k-1)+q$.
If  $k \geq 3$ is odd, then
 \[
 \ex(n,\h)=p\ex(k-1,K_r)+\ex(q,K_r).
 \]
If $k \geq 4$ is even, then
  \[
  \ex(n,\h) = \max \left\{p\ex(k-1,K_r)+\ex(q,K_r), \binom{\lfloor k/2 \rfloor -1}{2}+(\lfloor k/2 \rfloor -1)(n-\lfloor k/2 \rfloor +1)\right\}.
  \]
\end{theorem}
Apparently, for $k$ being odd and $k-1 | n$, Theorem \ref{conj11} gives that $\ex(n,\h)=\tfrac{n}{k-1}\ex(k-1,K_r)$ for all $n$. For $k$ being even and $k-1 | n$, one can verify
$$
\frac{n}{k-1}\ex(k-1,K_r)< \binom{\lfloor k/2 \rfloor -1}{2}+(\lfloor k/2 \rfloor -1)(n-\lfloor k/2 \rfloor +1)
$$
 for sufficiently large $n$ and then $\ex(n,\h)=\binom{\lfloor k/2 \rfloor -1}{2}+(\lfloor k/2 \rfloor -1)(n-\lfloor k/2 \rfloor +1)$ in this case. Therefore, we prove Conjecture \ref{conj1} in a stronger form as Theorem \ref{conj11} gives us the value of  $\ex(n,\h)$ for all $n$. In a recent paper, Fang, Zhu, and Chen \cite{FZC} verified  the exact version of Conjecture \ref{conj1} independently.

In order to prove Theorem \ref{thm1} and Theorem \ref{thm2}, we recall the study of $\ex(n, \C_{\geq k})$.
Actually, to determine the Tur\'an number of $\C_{\geq k}$, essentially one need  to consider the 2-connected analogue. To be precisely, let $\ex_{2\textrm{-conn}}(n,\C_{\geq k})$ be the maximum possible number of edges in an $n$-vertex 2-connected graph which is $\C_{\geq k}$-free.

Following the notation in the literature, for $2 \leq a \leq \lfloor (k-1)/2 \rfloor$, let $H(n,a,k)$ be a graph with vertex set $A \cup B \cup C$ such that $|A|=k-2a$, $|B|=a$, and $|C|=n-k+a$. Moreover, $A\cup B$ is a clique and each vertex in $B$ is adjacent to each vertex from $C$.
Let  $e(H(n,a,k))=h(n,a,k)=\binom{k-a}{2}+a(n-(k-a))$.  In 1976, Woodall \cite{W76} conjectured that $\ex_{2\textrm{-conn}}(n,\C_{\geq k})=\max\{h(n,2,k),h(n,\lfloor(k-1)/2 \rfloor,k)\}$ and verified this conjecture for large $n$. Fan, Lv, and Wang \cite{FLW04} proved the remaining case. Note that Kopylov  resolved this conjecture completely in the paper \cite{K77}.

For $\h=\{K_r, \C_{\geq k}\}$, let $\ex_{2\textrm{-conn}}(n,\h)$ be the maximum possible number of edges in a 2-connected $n$-vertex $\h$-free graph. For $2 \leq a \leq \lfloor (k-1)/2 \rfloor$, we vary the graph $H(n,a,k)$ through replacing the clique $K_{k-a}$ by the Tur\'an graph $T(k-a,r-1)$ and the resulting graph is denoted by $G_r(n,a,k)$. Let
$$g_r(n,a,k)=(n-k+a)a+\ex(k-a,K_r). $$
Note that  $g_r(n,a,k)=e(G_r(n,a,k))$ and $\ex(k-a,K_r)=\binom{k-a}{2}$ for  $ r>k-a$. We will prove the following results for  $\ex_{2\textrm{-conn}}(n,\h)$.
\begin{theorem} \label{thm11}
Let $\h=\{K_r, \C_{\geq k}\}$. If 	$\lfloor (k-1)/2 \rfloor+2 \leq r<k$ and $n \ge k\ge 5$, then
\[
\ex_{2\textrm{-conn}}(n,\h)=\max\left\{g_r(n,2,k),g_r(n,\lfloor (k-1)/2\rfloor,k) \right\}.
\]
\end{theorem}

For $r \leq \lfloor (k-1)/2 \rfloor+1$, our result is as follows.
	\begin{theorem} \label{thm21}
 Let $\h=\{K_r, \C_{\geq k}\}$.  If $r \le \lfloor (k-1)/2\rfloor+1$ and $n\ge k\geq 5$, then
\[
\ex_{2\textrm{-conn}}(n,\h)=e(G_2)
\]
 provided  $n \geq \tfrac{k^2}{2}$.
\end{theorem}
Here if $k=4$, then $r=2$ and the problem is not interesting. Thus we assume $k \geq 5$.

Since proofs of results in this paper involve different techniques, we will present them one by one so that each of them will be self-contained. The organization of  the rest of this paper as follows. In Section 2, we will prove Theorem \ref{thm11} and Theorem \ref{thm21}.
Proofs  of  Theorem \ref{thm1}, Theorem \ref{thm2}, and Theorem \ref{conj11}  will be presented in Section 3. In Section 4, we will  mention a few concluding remarks.
\section{Proofs of Theorem \ref{thm11} and Theorem \ref{thm21}}
\subsection{Proof of Theorem \ref{thm11}}
The proof of Theorem \ref{thm11} leverages on a variant of the method introduced by Kopylov \cite{K77}. Actually, the method from \cite{K77} turns out to be very useful in the study of Tur\'an type problems for graphs with bounded circumference, for example \cite{DNPWY,LYZ,L17,NP,zhang}.

We outline the proof of Theorem \ref{thm11} as follows. Assume $G$ is an edge-maximal counterexample and  $t=\lfloor (k-1)/2\rfloor$. Let $H(G,t)$ be
the $(t+1)$-core of $G$. We can show that $H(G,t)$ is not empty. Furthermore, we prove that $H(G,t)$ either is a clique or is a $K_r$-saturated graph. For the former case, we utilize a variant of the method  from \cite{K77}  to show the existence of a cycle with length at least $k$, which is a contradiction to the assumption of $G$. For the latter case, the key new idea is that we can treat $H(G,t)$ as a clique (see Lemma \ref{saturated}) and show the contradiction by repeating the argument in the former case.

  We  first prove a number of lemmas and the proof of Theorem \ref{thm11} will be presented at the end of Subsection 2.1. Readers can skip proofs of lemmas and go to the proof of Theorem \ref{thm11} directly.

 For the rest of this subsection, we fix
$$t=\lfloor (k-1)/2 \rfloor.$$
Recall the definition $g_r(n,k,a)=(n-k+a)a+\ex(k-a,K_r).$  We next show that $g_r$ is a convex function with $a$ being the independent variable.
\begin{lemma}\label{convex}
	The function $g_r(n,a,k)$ is  convex with $a$ being the independent variable, provided $r \geq 3$ and  $2 \leq a \leq t$.
\end{lemma}
\noindent
{\bf Proof:} Let $m(a)=g_r(n,a+1,k)-g_r(n,a,k)$. Then
\begin{align*}
	m(a)&=g_r(n,a+1,k)-g_r(n,a,k)\\
	&=(n-k+a+1)(a+1)+\ex(k-a-1,K_r)-(n-k+a)a-\ex(k-a,K_r)\\
	&=n-k+2a+1+\big(\ex(k-a-1,K_r)-\ex(k-a,K_r)\big).
\end{align*}
To get $T(k-a,r-1)$, one can add a new vertex to a smallest partite of $T(k-a-1,r-1)$ and connect this vertex to all vertices from other parties. Here the smallest partite is an empty set if $k-a-1<r-1$. Thus
$$
\ex(k-a,K_r)-\ex(k-a-1,K_r)=k-a-1-\lfloor(k-a-1)/(r-1) \rfloor.
$$
It follows that
\begin{align*}
	m(a+1)-m(a)&=n-k+2(a+1)+1+\big(\ex(k-a-2,K_r)-\ex(k-a-1,K_r)\big)\\
	&\quad -(n-k+2a+1)-\big(\ex(k-a-1,K_r)-\ex(k-a,K_r)\big)\\
	&=2-(k-a-2-\lfloor (k-a-2)/(r-1)\rfloor)+(k-a-1-\lfloor(k-a-1)/(r-1)\rfloor )\\
	&=3-\lfloor(k-a-1)/(r-1)\rfloor+ \lfloor (k-a-2)/(r-1)\rfloor\\
	&\ge 3-(k-a-1)/(r-1)+(k-a-2)/(r-1)-1\\
	&=2-\frac{1}{r-1}>0.
\end{align*}
Thus $g_r(n,k,a)$ is convex if $r \geq 3$ and $2 \leq a \leq t$.  \hfill $\square$

We will need the following two results on the length of a longest cycle in a graph.
\begin{theorem}[Dirac \cite{D52}] \label{dirac}
Let $G$ be a 2-connected graph of order $n$ with minimum degree $\delta$. Then $c(G)\geq \mbox{min}\{2\delta,n\}$, where $c(G)$ is the length of a longest cycle in $G$.
\end{theorem}

\begin{lemma}[Kopylov \cite{K77}]\label{klemma}
	Let $G$ be a $2$-connected $n$-vertex graph with a path $P$ of $m$ edges with endpoints $x$ and $y$. For $v\in V(G)$, let $d_P(v)=|N(v)\cap V(P)|$. Then $G$ contains a cycle of length at least $\min\{m+1,d_P(x)+d_P(y)\}.$
\end{lemma}

A graph $G$ is $K_r$-{\it saturated} if $G$ is $K_r$-free and adding any nonedge to $G$ will create a copy of $K_r$. In other words,  any two nonadjacent vertices are in a $K_r-e$, here $K_r-e$ is the clique $K_r$ with one edge removed. The next lemma plays an important role in our proof of Theorem \ref{thm11}.
\begin{lemma}\label{saturated}
Assume that  $H$ is an induced subgraph of $G$  and $H$ is  $K_r$-saturated with $\delta(H) \geq \lfloor (k-1)/2 \rfloor+1$ and $k \geq 5$. Let $G'$ be a graph by adding all nonedges with endpoints in $H$. If $r \geq \lfloor (k-1)/2 \rfloor+2$ and the length of a longest cycle in $G'$ is  at least $k$, then there is a cycle of length at least $k$ in $G$.
\end{lemma}
\noindent
{\bf Proof:}  Let $\overline{H}$  be the complement graph of $H$. Among all cycles of length at least $k$ in $G'$, we pick $C$ as the one satisfies the following two conditions:

{ (a):} $|E(C)\cap E(\overline{H})|$ is minimized; and

{ (b):} the length is maximized, subject to (a).

Let $C=a_1a_2\cdots a_\ell a_1$ with $\ell \geq k$ and vertices of $C$ are ordered clockwisely.
For each $1 \leq i \leq \ell$, we will write $a_i^+$ for $a_{i+1}$ and $a_i^-$ for $a_{i-1}$ in some circumstances, here the addition is under modulo $\ell$. We claim that $C$ is a cycle in $G$.  If $C$ contains edges from $E(\overline{H})$, then we next find a contradiction to the definition of $C$. By the construction of $G'$, we may assume that $a_1a_{\ell} \in E(\overline{H})$. Note that $\{a_1,a_\ell\} \subset V(H)$. Define  $N_H(a)=N_G(a)\cap V(H)$.  Assume that $N_C(a_1) \cap N_C(a_{\ell})=\{a_{i_1},\ldots,a_{i_q}\}$ with vertices ordered clockwisely. Let $P_1=a_{i_1} a_{i_1}^+ \cdots a_{i_q}$, $P_2=a_1a_2 \cdots a_{i_q}$, and  $P_3=a_{i_1} a_{i_1}^+ \cdots a_{\ell}$.

 \begin{claim} \label{claim1}
 $(1)$  $N_G(a_1) \cap N_G(a_{\ell}) \subset V(C)$, $N_H(a_1) \subset V(C)$, and $N_H(a_{\ell}) \subset V(C)$.\\
$(2)$  For each $1 \leq j \leq q-1$, $a_{i_j}^+ \neq a_{i_{j+1}}$. \\
$(3)$ $a_1 a_{i_j}^+ \not \in E(G)$ for each $1 \leq j \leq q$ and $a_\ell a_{i_j}^- \not \in E(G)$ for each $1 \leq j \leq q$.\\
$(4)$ $|V(P_1)|\in\{k-3,k-2\},|V(P_2)|\le k-1$, and $|V(P_3)|\le k-1.$
\end{claim}
\noindent
{\bf Proof of Claim \ref{claim1}:} For part $(1)$,  if there is a vertex $v\in N_G(a_1)\cap N_G(a_\ell)$ such that  $v \not \in V(C)$, then we  get a new cycle $C_1= va_1\cdots a_{\ell-1}a_{\ell}v$ of length at least $k$  and $|E(C_1)\cap E(\overline{H})|<|E(C)\cap E(\overline{H})|$. This is a contradiction to the choice of $C$.  Therefore, $N_G(a_1)\cap N_G(a_{\ell}) \subset V(C)$. Suppose that there is a vertex $v\in N_H(a_1)$ such that  $v \not \in V(C)$.
As $H$ is a clique in $G'$ and $v, a_\ell \in V(H)$,  then $va_1 \cdots a_\ell v$ is a cycle with at least $k$ edges in $G'$ .  If $va_\ell \in E(G)$, then it violates the condition (a).   If $va_\ell \not \in E(G)$, then it violates the condition (b). Thus  $N_H(a_1) \subset V(C)$. By symmetry, we also have $N_H(a_{\ell}) \subset V(C)$.

For part $(2)$,  if  there is a $a_{i_j}$ such that $a_{i_j}^+ = a_{i_{j+1}}$, then we get a new cycle $C_1=a_1a_2\cdots a_{i_j}a_\ell a_\ell^-\cdots  a_{i_{j+1}}a_1$ with at least $k$ edges and  $|E(C_1)\cap E(\overline{H})|<|E(C)\cap E(\overline{H})|$. This is a contradiction to the definition of $C$.

The proof of part $(3)$ is the same as the one for part $(2)$ and it is skipped here.

For part $(4)$, we first claim $|V(P_1)|\in\{k-3,k-2\}$. If $|V(P_1)|\ge k-1$,  then $C_1=a_1 a_2 \ldots a_{i_q}a_1$ is a new cycle with at least $k$ edges and $|E(C_1)\cap E(\overline{H})|<|E(C)\cap E(\overline{H})|$, here note that $a_1a_{i_q} \in E(G)$. This is a contradiction to the definition of $C$. Next we claim $|V(P_1)|\ge k-3.$  Recall  $N_C(a_1)\cap N_C(a_{\ell})=\{a_{i_1},\ldots,a_{i_q}\}$. By part $(1)$, we have $N_G(a_1)\cap N_G(a_{\ell}) \subset V(C)$.
As $H$ is $K_r$-saturated, any two nonadjacent vertices in $H$ is contained in a $K_r-e$. Notice that$\{a_1,a_{\ell}\} \subset V(H)$ and $a_1$ is not adjacent to $a_{\ell}$. Thus
\[
q=|N_G(a_1)\cap N_G(a_{\ell})|\ge r-2.
\]
By part $(2)$, we have  $a_{i_j}^+ \neq a_{i_{j+1}}$, i.e., $i_{j+1}-i_j \geq 2$. One can observe that
\[
V(P_1) \geq 2(q-1)+1.
\]
Recall that $q \geq r-2 \geq \lfloor (k-1)/2 \rfloor$. It follows that $V(P_1) \geq 2\lfloor (k-1)/2 \rfloor-1 \geq k-3.$ Thus $|V(P_1)|\in\{k-3,k-2\}.$
Finally, we claim $|V(P_2)|\le k-1$ and $|V(P_3)|\le k-1$. By symmetry, one can suppose that $|V(P_2)| \geq k$. Now $C_1=a_1a_2\cdots a_{i_q}a_1$ is a new cycle with at least $k$ edges and $|E(C_1)\cap E(\overline{H})|<|E(C)\cap E(\overline{H})|$. This is a contradiction to the definition of $C$.
\hfill$\square$

We continue to prove the lemma.
For the case of $r \geq \lfloor (k-1)/2 \rfloor+3$, it follows that  $q \geq r-2 \geq \lfloor (k-1)/2 \rfloor+1$ and $V(P_1) \geq 2\lfloor (k-1)/2 \rfloor+1 \geq k-1.$ This is a contradiction to part $(4)$ of Claim 1.
It remains to show the desired contradiction for the case where $r=\lfloor (k-1)/2 \rfloor+2$. Repeating the argument above, we can observe that $q=\lfloor (k-1)/2 \rfloor$. Notice that  $|V(P_1)|\in\{k-3,k-2\}$ by part $(4)$ of Claim 1.

If $|V(P_1)|=k-2$, then $\ell=k$, $a_{i_1}=a_2$, and $a_{i_q}=a_{k-1}$ by $(4)$ of Claim 1. Recall $|V(P_1)|\ge 2(q-1)+1$. For the case that  $k$ is odd, then $|V(P_1)|\ge k-2$. This implies $a_{i_{j+1}}=a_{i_j}+2$ for any $1\le i\le q-1.$ Since $|N_H(a_1)| \geq \lfloor (k-1)/2 \rfloor+1$,  $N_H(a_1) \subset V(C)$ by part $(1)$ of Claim 1, and $|N_G(a_1)\cap N_G(a_\ell)|=q=\lfloor (k-1)/2 \rfloor$,  it follows that there is some  $a_{i}$ such that $a_i \in N_H(a_1)$ but $a_i \not \in N_G(a_1)\cap N_G(a_\ell)$. The only possibility is that $a_i=a_{i_j}^+$ for some $a_{i_j}$. However, this is a contradiction to part  $(3)$ of  Claim 1. For the case that $k$ is even, there exists a unique $1 \leq j \leq q$ such that  $a_{i_{j+1}}=a_{i_j}+3$. Similarly, there are $a_i$ and $a_j$ such that $a_i \in N_H(a_1) \setminus N_G(a_{\ell})$ and
 $a_j \in N_H(a_\ell) \setminus N_G(a_{1})$. Part $(3)$ of Claim 1 yields that the only possibility is $a_i=a_{i_{j+1}}^-$ and $a_j=a_{i_j}^+$. Now we can get a new cycle $C_1=a_1a_{i_{j+1}}^- a_{i_{j+1}} \cdots a_\ell a_{i_j}^+ a_{i_j} \cdots a_1$ with at least $k$ edges and $|E(C_1)\cap E(\overline{H})|<|E(C)\cap E(\overline{H})|$. This is a contradiction to the definition of $C$.

  For $|V(P_1)|=k-3$, we get that $k$ is even and   $a_{i_{j+1}}=a_{i_j}+2$ for any $1\le i \le q-1$.
 Recall part $(4)$ of Claim 1 and the assumption $|V(C)| \geq k$.  Either $\ell=k,a_{i_1}=a_3, a_{i_q}=k-1$, or $\ell=k,a_{i_1}=a_2,a_{i_q}=k-2$, or $\ell=k+1,a_{i_1}=a_3,a_{i_q}=k-1$, or $\ell=k+1,a_{i_1}=a_2,a_{i_q}=k-2.$ By symmetry, we need only to consider the first case and the third case.  For the first case,  there is a vertex $a_i \in N_H(a_{\ell}) \setminus N_G(a_1)$ as the minimum degree assumption. It is only possible that $a_i=a_{i_j}^-$ for some $a_{i_j}$. This is a contradiction to
part $(2)$ of Claim 1.  For the third case, consider the cycle $C_1=a_1 \cdots a_{i_1} a_{\ell} \cdots a_{i_2} a_1$. Note that $C_1$ contains $k$ edges and $|E(C_1)\cap E(\overline{H})|<|E(C)\cap E(\overline{H})|$. This is a contradiction to the definition of $C$.   The lemma is proved. \hfill $\square$

We are ready to prove Theorem \ref{thm11}.

\vspace{0.2cm}
\noindent
{\bf Proof of Theorem \ref{thm11}}: Notice that $r \geq t+2$.  For the lower bound,
we recall  the definition of $G_r(n,a,k)$. It is obvious that $G_r(n,a,k)$ is $\{K_r,\C_{\geq k}\}$-free and the lower bound follows.
We mention a special case where $k$ is even, $r=t+2$,  and $a=t=k/2-1$.
In this case, the subgraph of $G_r(n,a,k)$ induced by  $A \cup B$ is $K_r$ with one edge removed.

 For the upper bound, if the assertion in Theorem \ref{thm11} is not true, then let $G$ be an edge-maximal counterexample. That is $G$ is $\h$-free and
$$
e(G)> \max\{g_r(n,2,k),g_r(n,t,k)\}.
$$
 Moreover, adding any nonedge to $G$ either creates a $K_r$ or results a cycle with length at least $k$.

Let $H(G,t)$ be the $(t+1)$-core of $G$, i.e., the largest induced subgraph of $G$ with minimum degree at least $t+1$.

\begin{claim}
	$H(G,t)$ is not empty.
\end{claim}
\noindent
{\bf Proof:} Note that if $k$ is odd, then the subgraph of $G_r(n,t,k)$ induced by $A\cup B$ is a clique of size $t+1$ and each vertex in $C$ is adjacent to all vertices in $B$. If $k$ is even, then the subgraph induced by $A \cup B$ contains $K_{t+1}$ with one edge removed as a subgraph. Therefore, $e(G_r(n,t,k))=g_r(n,t,k) \geq (n-t)t+\binom{t}{2}$ and $e(G)>(n-t)t+\binom{t}{2}$ by the assumption. If $H(G,t)$ is empty, then in the process of defining $H(G,t)$, we remove  at most $t$ edges for each of the first $n-t$ vertices  and the subgraph of the last $t$ vertices contains at most $\binom{t}{2}$ edges. Thus $e(G) \leq t(n-t)+\binom{t}{2}$ and this is  a contradiction to the lower bound for $e(G)$. \hfill $\square$

\begin{claim}
	If $x,y \in H(G,t)$ and $x$ is not adjacent to $y$, then there is no $xy$-path with at least $k-1$ edges in $G$.
\end{claim}
\noindent
{\bf Proof:} If the assertion does not hold, then we choose vertices $x$ and $y$ from $H(G,t)$ such that a longest path in $G$ from
$x$ to $y$ contains the largest number of edges among all such nonadjacent pairs. Let $P$ be a longest $xy$-path in $G$.
We assert that all neighbors of $x$ in $H(G,t)$ lie in $P$. If $x$ has a neighbor $x' \in H(G,t)\setminus P$, then  either $x'y \in E(G)$ or $x'y \notin E(G)$. In the former case, $x'xPyx'$ is a cycle of length at least $k$ in $G$, a contradiction to the assumption for $G$. For the latter case, $x'xPy$ is a longer path, which is a contradiction to the selection of $x$ and $y$. There is a contradiction in each case. We can show the same assertion for $y$ similarly. Therefore, $d_P(x) \geq t+1$ and $d_P(y) \geq t+1$. By Lemma \ref{klemma},  $G$ contains a cycle of length at least $\min\{k,d_P(x)+d_P(y)\}=k$, a contradiction. \hfill $\square$

The next claim follows from the assumption for $G$.
\begin{claim}
	$H(G,t)$ is either a clique or a $K_r$-saturated graph.
\end{claim}
\begin{claim}
	$t+2\le h=|H(G,t)|\le k-2.$
\end{claim}
\noindent
{\bf Proof:} The lower bound is trivial as the minimum degree assumption. Note that $r \geq 4$ because we assume $k \geq 5$ and $r \geq t+2$. As $H(G,t)$ is either a clique or  $K_r$-saturated, then it is 2-connected. If $h \geq k$, then Theorem \ref{dirac} ensures that $H(G,t)$ contains a cycle of length at least $\min\{h,2(t+1)\} \geq k$, a contradiction. If $h=k-1$, then a classical result by Ore \cite{ore} implies that $H(G,t)$ is Hamiltonian connected . Since $H(G,t)$ is an induced subgraph of a 2-connected graph $G$, for any vertex $x$ not in $H(G,t)$, there are two internally disjoint paths from $x$ to $H(G,t)$. Let $y$ and $z$ be the endpoints of these two paths in $H(G,t)$. As there is a  Hamilton path from $y$ to $z$ in $H(G,t)$, we get that $G$ contains a cycle with at least $k$ edges, a contradiction.  \hfill $\square$

Let $H(G,k-h)$ be the $(k-h+1)$-core of $G$.
\begin{claim}
	$H(G,t)\neq  H(G,k-h)$.
\end{claim}
\noindent
{\bf Proof:} Observe that $2 \leq k-h \leq t$. If 	$H(G,t)=H(G,k-h)$, then in the process of defining $H(G,t)$, we remove at most $k-h$ edges for the first $n-h$ vertices and the subgraph in the last $h$ vertices contains at most $\ex(h,K_r)$ edges as $G$ is $K_r$-free. We get that
\[
e(G) \leq \ex(h,K_r)+(n-h)(k-h)=g_r(n,k-h,k) \leq \max\{g_r(n,2,k),g_r(n,t,k)\},
\]
a contradiction. Here the last inequality holds as the range of $k-h$ and Lemma \ref{convex}. \hfill $\square$

We first consider the case where $H(G,t)$ is not a clique, i.e., $H(G,t)$ is $K_r$-saturated.  Denote $X=V(H(G,t))$ and $Y=V(H(G,k-h))$. As we already proved $H(G,t) \neq H(G,k-h)$, it follows that $X$ is a proper subset of  $Y$. In addition, for each $y \in Y \setminus X$, there is  $x \in X$ such that $x$ and $y$ are not adjacent in $G$. Otherwise, $d_X(y) \geq h \geq t+2$ and $y$ should be a vertex in $X$, a contradiction.

If  $r \geq t+3$, then for any nonadjacent $x$ and $y$ with $x \in X$ and $y \in Y \setminus X$,
we claim that adding the edge $xy$ to $G$ must create a cycle with length at least $k$. Otherwise, $x$ and $y$ are contained in $K_r-e$ by the assumption of $G$. As $r \geq t+3$, all vertices in $K_r-e$ has at least $t+1$ neighbors  in $K_r-e$. Then $y \in H(G,t)$, a contradiction. Let
\[
Q=\{(x,y): x \in X, y \in Y\setminus X, \textrm{ and } x \textrm{ is not adjacent to } y\}.
\]
Let $G'$ be a new graph by adding all nonedges in $X$. For each pair $(x,y) \in Q$, we already showed that there is an $xy$-path  with at least $k-1$ edges in $G$. Then the same observation also holds in $G'$. We choose the pair $(x,y) \in Q$ such that a longest path $xy$-path, say $P_{xy}$,  in $G'$  contains the largest number of edges among all pairs in $Q$.
Notice that $N_{G'}(y)=N_G(y)$. We claim that $N_{G'}(y) \cap Y \subset V(P_{xy})$. If not, then let $z \in N_{G'}(y) \cap Y $ and $z \not \in V(P_{xy})$. As the choice of the pair $(x,y)$, we get that $x$ and $z$ must be adjacent. Thus $xP_{xy}yzx$ is a cycle of length at least $k$ in $G'$. By Lemma \ref{saturated}, $G$ contains a cycle of length $k$, a contradiction. Similarly, we can show $X \subset V(P_{xy})$. Therefore, by Lemma \ref{klemma}, $G'$ contains a cycle of length at least $\min\{k, h-1+k-h+1\}=k$, here note that $X$ is a clique in  $G'$. Then $G$ contains a cycle with at least $k$ edges by Lemma \ref{saturated}. We obtained a contradiction in each case.

For $r=t+2$, one can observe that each vertex  $y \in Y\setminus X$ has at most $t$ neighbors in $X$. Otherwise, $y \in X$. Claim 5 implies that each vertex from $Y\setminus X$ has at least two nonneighbors in $X$. Furthermore, if there are vertices  $x \in X$ and $y \in Y \setminus X$ such that $xy \not \in E(G)$ and
adding the edge $xy$ to $G$  creates a cycle with length at least $k$, then one can obtain the same contradiction by reusing the argument above. Therefore, for any  nonadjacent vertices $x$ and $y$ such that $x \in X$ and $y \in Y \setminus X$,  adding the edge $xy$ to $G$ must create a copy of $K_{t+2}$, where the vertex set of $K_{t+2}$ is $\{x,y\} \cup T$ for some $t$-set $T$.  Notice that $K_{t+2}-xy$ is a subgraph induced by $Y$.  It follows that $T \subset Y$ and  $d_Y(y) \geq t$.

 We next show $Y \setminus X$ is an independent set. It only needs to consider the case of $|Y \setminus X| \geq 2$. Suppose that there is an edge $y_1y_2$ with $\{y_1,y_2\} \subset Y \setminus X$. For $i=1,2$, let $x_i \in X$ be a vertex such that $x_i$ and $y_i$ are not adjacent. Additionally, by the assumption, there is a $t$-set $T_i \subset Y$ such that $\{x_i,y_i\} \cup T_i$ induces a $K_{t+2}-e$. We claim that either $y_2 \in T_1$ or $y_1 \in T_2$. Otherwise,
 observe that the subgraph induced by $\cup_{i=1}^2 \{y_i \cup T_i\} \cup X$ has minimum degree at least $t+1$. Thus $y_i \in X$.  This is a contradiction to the choice of $y_i$. Without loss of generality, we assume $y_2 \in T_1$.   One can see $y_1$ is not adjacent to any vertex from $X \setminus T_1$. Otherwise,  the subgraph induced by $y_1\cup T_1 \cup X$ has minimum degree at least $t+1$. It implies that $y_1 \in X$, which is  a contradiction. For any $x_1' \in X\setminus T_1$ which is a nonneighbor of $y_1$,  let $T_1'$ be  a $t$-set from $Y$ such that $\{y_1,x_1'\} \cup T_1'$ is a $K_{t+2}-e$. Then $T_1=T_1'$. Otherwise, the subgraph induced by $y_1 \cup T_1 \cup T_1' \cup X$ has minimum degree at least $t+1$ and then $y_1 \in X$, a contradiction. Therefore, any vertex $x_1 \in X\setminus T_1$ together with $y_1 \cup T_1$ induce a $K_{t+2}-e$. It follows that $X \setminus T_1$ is an independent set as $G$ is $K_{t+2}$-free. Equivalently, $N_X(x_1) \subset X \cap T_1$, which yields that $d_X(x_1) \leq t$ for each $x_1 \in X\setminus T_1$. This is a contradiction to the definition of $X=H(G,t)$. Thus $Y \setminus X$ is an independent set.

 Fix nonadjacent vertices $y \in Y\setminus X$ and $x \in X$.
Let $T \subset Y$ be the vertex of a $K_t$ such that $\{x,y\} \cup T$ induces $K_{t+2}-e$. As $Y \setminus X$ is an independent set, we get that $T \subset X$. Moreover, $N_Y(y)=T$ and $y$ is not adjacent to each $x' \in X\setminus T$. Otherwise, $y$ has at least $t+1$ neighbors in $X$ and $y$ should be a vertex in $X$, a contradiction.
 Therefore, $\{x',y\} \cup T$ is $K_{t+2}-e$ for any $x' \in X \setminus T$. Since $G$ is $K_{t+2}$-free, $X \setminus T$ is an independent set. This implies that $d_{X}(x')=t$ for each $x' \in X \setminus T$, a contradiction to the definition of $H(G,t)$.

Next, we consider the case  where $H(G,t)$ is a clique. If $r \geq t+3$, then we repeat the argument above by defining $G'=G$ and get a contradiction. For  $r=t+2$, since $t+2\le |H(G,t)|\le k-2$, then $H(G,t)$ contains a $K_r$, a contradiction. The proof of Theorem \ref{thm11} is complete. \hfill $\square$ 

\subsection{Proof of Theorem \ref{thm21}}
As we are not able to prove Lemma \ref{saturated} under assumptions in Theorem \ref{thm21}, we apply a different approach to
 prove Theorem \ref{thm21}. Namely, let $G$ be a 2-connected graph which is both $K_r$-free and $\C_{\geq k}$-free.  A longest cycle in $G$ is denoted by $C$. Let $e(C)$ be the number of edges contained in $C$ and $e(G-C)+e(G-C,C)$ be the number of edges which has at most one endpoint in $C$. Apparently,
 $$e(G)= e(C)+e(G-C)+e(G-C,C).$$
 We will upper bound $e(C)$ and $e(G-C)+e(G-C,C)$ separately.  Lemma \ref{newstability} will give us an upper bound for $e(G)$.

We start with the following estimate for the number of edges in $C$.
\begin{lemma}\label{cycle edge}
	Let $G$ be a $K_r$-free graph and $C$ be a longest cycle in $G$ of length $c$, where $r\le \lfloor c/2 \rfloor+1$.  If there exists a vertex $u\in V(G-C)$ with $d_C(u)=\lfloor c/2 \rfloor$, then
	$$e(C)\le\ex(\lfloor c/2 \rfloor,K_{r-1})+\lfloor c/2 \rfloor\lceil c/2 \rceil.$$
\end{lemma}
\noindent
{\bf Proof:}  If $c$ is even, then let  $C=a_1b_1a_2b_2\cdots a_{c/2}b_{c/2}a_1$ be the longest cycle. Since $d_C(u)= c/2$, without loss of generality, assume $A=V(C)\cap N_C(u)=\{a_1,\ldots,a_i,\ldots,a_{c/2}\}$. Let $B=\{b_1,\ldots,b_i,\ldots,b_{c/2}\}$. Then $e(B)=0$. Otherwise, there is a longer cycle by including $u$, a contradiction.  As $G[C]$ is also $K_r$-free, we get $e(A) \leq \ex(c/2,K_{r-1})$. Thus
\begin{align*}
	e(C)&=e(A)+e(B)+e(A,B)\\
	&\le \ex(c/2,K_{r-1})+(c/2)^2 \\
    &=\ex( c/2 ,K_{r-1})+ (c/2)^2.
\end{align*}
If $c$ is odd, then let   $C=a_1b_1a_2b_2\cdots a_{\lfloor c/2\rfloor}b_{\lfloor c/2\rfloor}b_{\lceil c/2\rceil}a_1$ be the longest cycle. Similarly, we assume $A=V(C)\cap N_C(u)=\{a_1,\ldots,a_i,\ldots,a_{\lfloor c/2\rfloor}\}$ and $B=\{b_1,\ldots,b_i,\ldots,b_{\lfloor c/2\rfloor},b_{\lceil c/2\rceil}\}$. Then $G[B]$  contains at most one edge, that is $b_{\lfloor c/2\rfloor}b_{\lceil c/2\rceil}$. Otherwise, there is a longer cycle by including $u$, a contradiction.  Thus
\begin{align*}
	e(C)&=e(A)+e(B)+e(A,B)\\
	&\le\ex(\lfloor c/2\rfloor,K_{r-1})+1+\lfloor c/2 \rfloor\lceil c/2 \rceil.
\end{align*}
The inequality above must be strict as $r \le\lfloor c/2 \rfloor+1$. Otherwise, a $K_{r-2}$ in $A$ together with vertices $b_{\lfloor c/2\rfloor}$ and $b_{\lceil c/2\rceil}$ form a $K_r$, a contradiction. Thus
\[
e(C)\le\ex(\lfloor c/2\rfloor,K_{r-1})+\lfloor c/2 \rfloor\lceil c/2 \rceil
\]
in this case.
 The lemma is proved. \hfill$\square$

Let us recall the following classical result by Bondy.
\begin{lemma}[Bondy \cite{Bondy}]\label{Bondy}
	Let $G$ be a graph on $n$ vertices and $C$ be a longest cycle of $G$ with order $c$. Then
	\[
	e(G-C)+e(G-C,C)\le\lfloor c/2 \rfloor(n-c).
	\]
\end{lemma}
 Woodall \cite{W76} made a conjecture on the maximum number of edges in a 2-connected graph which is $\C_{\geq k}$-free and has minimum degree at least $k \geq 2$. Ma and Ning \cite{MN} proved a stability result of Woodall's conjecture. Previously, F\"uredi, Kostochka, and Verstra\"ete \cite{FKV} obtained a stability result for  Erd\H{o}s-Gallai theorems on paths and cycles (see Theorem \ref{EGpath} and Theorem \ref{EGcycle}). In \cite{MN}, Ma and Ning discovered  the stability of Woodall's conjecture is directly related to the above classical result of Bondy. Indeed, they proved a stability result of Lemma \ref{Bondy}, see Theorem 3.2 in \cite{MN}. For our problem, we need the following variant of Theorem 3.2 in \cite{MN}.
\begin{lemma}\label{stability}
	Let $G$ be a graph on $n$ vertices and $C$ be a longest cycle in $G$ of length $c$, where $4\le c\le n-1$. If $e(G-C)+e(G-C,C)>(\lfloor c/2 \rfloor-\tfrac{1}{2})(n-c)$, then one of the following holds:
	
	\noindent
	$(1)$ there exists a vertex $u\in V(G-C)$ with $d_C(u)=\lfloor c/2 \rfloor;$\\
	\noindent
	$(2)$ there exists a cycle $C'$in $G$ such that $|V(C\cap C')\le 1$. Moreover, if $V(C\cap C')=\emptyset$, then $|C'|\ge 2\lfloor\frac{c}{2}\rfloor-3.$ If $|V(C\cap C')|=1$, then $|C'|\ge 2\lfloor\frac{c}{2}\rfloor-1.$
\end{lemma}
The proof of the  lemma above is similar to the one for Theorem $3.2$ in Ma and Ning's paper \cite{MN} and it is omitted here. If $G$ is 2-connected, then Ma and Ning pointed out that $(2)$  does not occur. Thus the combination of  Lemma \ref{Bondy} and Lemma \ref{stability} gives the following lemma.
\begin{lemma}\label{g-c}
	Let $G$ be a $2$-connected graph on $n$ vertices and $C$ be a longest cycle in $G$ of length $c$ with $4 \leq c\le n-1$.
	
	\noindent
	$(1)$ If there exists a vertex $u\in V(G-C)$ with $d_C(u)=\lfloor c/2 \rfloor$, then
	\[e(G-C)+e(G-C,C)\le\lfloor c/2 \rfloor(n-c).\]
	\noindent
	$(2)$ If $d_C(u)<\lfloor c/2 \rfloor$ for each  vertex $u\in V(G-C)$, then
	\[e(G-C)+e(G-C,C)\le(\lfloor c/2 \rfloor-\tfrac{1}{2})(n-c).\]
\end{lemma}

Let $t=\lfloor (k-1)/2\rfloor$. The following inequality will be used in the proof of Lemma \ref{newstability}.
\begin{lemma} \label{k-2}
		If $k\le n$ and $r\ge 3$, then
	\[
	\ex(k-2,K_r)+\left(t-\frac{1}{2}\right)(n-k+2) <\ex(k-1,K_r)+\left(t-\frac{1}{2}\right)(n-k+1).
	\]
\end{lemma}
{\bf Proof:}  Notice that $\ex(k-1,K_r)$  and $\ex(k-2,K_r)$ are the number of edges in the Tur\'an graph $T(k-1,r-1)$ and $T(k-2,K_r)$ respectively. The Tur\'an graph $T(k-2,r-1)$ can be obtained by deleting a  vertex  from the largest partite  of graph  $T(k-1,r-1)$. This implies $\ex(k-1,K_r)-\ex(k-2,K_r)=k-1-\lceil(k-1)/(r-1)\rceil$.   Thus it is sufficient to  show
\[
k-1-\lceil(k-1)/(r-1)\rceil > t-1/2= \lfloor(k-1)/2\rfloor-1/2.
\]
 As $r\ge 3$, then $k-1-\lceil(k-1)/(r-1)\rceil\ge k-1-\lceil(k-1)/2\rceil= \lfloor(k-1)/2\rfloor$. The  inequality above is true. \hfill$\square$

 The next lemma  gives us an upper bound for the number of edges in a 2-connected $\{K_r,\C_{\geq k}\}$-free graph which is of independent interest.
\begin{lemma}\label{newstability}
Assume that $n\ge k \geq 5$ and $r \leq t+1$. Let $G$ be a $2$-connected  $\{K_r,\C_{\geq k}\}$-free graph with $n$ vertices and $C$ be a longest cycle in $G$.

\noindent
$(1)$ If $|V(C)|=k-1$  and there exists $u \in V(G)\setminus V(C)$ such that $u$ has $t$ neighbors in $C$, then
	$$
	e(G) \leq\ex(t,K_{r-1})+t(n-t).
	$$
\noindent
$(2)$	If $k$ is even, $|V(C)|=k-2$, and there exists $u \in V(G)\setminus V(C)$ such that $u$ has $t$ neighbors in $C$, then
	$$
	e(G) \leq\ex(t,K_{r-1})+t(n-t).
	$$
\noindent
$(3)$ Otherwise, we have
	\[
	e(G) \leq \max\{(t-1/2) (n-1),\ex(k-1,K_r)+(t-1/2)(n-k+1)\}.
	\]
\end{lemma}
\noindent
{\bf Proof:}  Note that $G$ is a  2-connected $n$-vertex graph that is $\{K_r, \C_{\geq k}\}$-free. In addition, $C$ is a longest cycle in $G$. It implies that $4 \leq c \leq k-1 \leq n-1$. The lower bound for $c$ follows  from the assumption that $G$ is 2-connected.

We start to prove part $(1)$. Notice that $|V(C)|=k-1$ and   there is  a vertex $u\in V(G\setminus C)$ such that $d_C(u)=t$. By Lemma \ref{g-c}, we have
\[
e(G-C)+e(G-C,C)\le t(n-k+1).
\]
Lemma \ref{cycle edge} tells us that
\begin{align*}
	e(C)+e(G-C)+e(G-C,C)&\le\ex(t,K_{r-1})+t(k-1-t)+t(n-k+1)\\
	&=\ex(t,K_{r-1})+t(n-t).
\end{align*}
Thus the assertion in $(1)$ follows.

For part $(2)$, notice that $t=(k-2)/2$ as $k$ is even. Repeating the proof for part $(1)$, we can show part $(2)$.

For part $(3)$,  if  $|V(C)|=k-1$, then  each $u\in V(G\setminus C)$ has at most $t-1$ neighbors in $C$. By Lemma \ref{g-c}, we have
\[
e(G-C)+e(G-C,C)\le\left(t-\tfrac{1}{2}\right)(n-k+1).
\]
Apparently, $e(C) \leq \ex(k-1,K_r)$. Thus
\[
	e(C)+e(G-C)+e(G-C,C) \le\ex(k-1,K_r)+\left(t-\tfrac{1}{2}\right)(n-k+1).\\
\]

If $|V(C)| \leq k-3$, then   Theorem \ref{EGcycle} yields that
\[
e(G) \le \frac{k-3}{2}(n-1)  \leq (t-\frac{1}{2})(n-1)
\]
no matter the parity of $k$.  If $|V(C)|= k-2$ and $k$ is odd, then $e(G) \leq (t-\tfrac{1}{2})(n-1)$ also follows from  Theorem \ref{EGcycle}. If $|V(C)|= k-2$ and $k$ is even, then  each   vertex from $u\in V(G\setminus C)$ has at most $t-1$ neighbors in $C$. By Lemma \ref{g-c}, we get that
\[
e(G-C)+e(G-C,C)\le(t-\tfrac{1}{2})(n-k+2).
\]
Then $e(G)$ can be upper bounded as follows:
\begin{align*}
	e(G)&=e(C)+e(G-C)+e(G-C,C)\\
	&\le \ex(k-2,K_r)+(t-\tfrac{1}{2})(n-k+2)\\
	&< \ex(k-1,K_r)+(t-\frac{1}{2})(n-k+1),
\end{align*}
the last inequality is true by  Lemma \ref{k-2}. Part $(3)$ is proved. \hfill $\square$

We are now ready to prove Theorem \ref{thm21}.

\vspace{0.2cm}
\noindent
{\bf Proof of Theorem \ref{thm21}:} For the lower bound, recall that $t=\lfloor (k-1)/2\rfloor$ and the graph $G_2$.
Observe that the clique number of $G_2$ is $r-1$ and a longest cycle has $2t \leq k-1$ edges. Moreover, this graph is 2-connected. The lower bound is proved.

For the upper bound, let $G$ be an $n$-vertex 2-connected graph which is $\{K_r, \C_{\geq k}\}$-free.
By Lemma \ref{newstability}, either
\[
e(G) \leq \ex(t,K_{r-1})+t(n-t),
\]
or
\[
e(G) \leq \max\{(t-1/2) (n-1),\ex(k-1,K_r)+(t-1/2)(n-k+1)\}.
\]
We need only to consider the latter case.  Observe that
\begin{align*}
e(G) &\leq \max\{(t-1/2) (n-1),\ex(k-1,K_r)+(t-1/2)(n-k+1)\}\\
&\le  \ex(t,K_{r-1})+t(n-t)
\end{align*}
as $n\ge \tfrac{k^2}{2}$.  The proof for Theorem \ref{thm21} is complete. \hfill$\square$

%

%

%
\section{Proofs of Theorem \ref{thm1}, Theorem \ref{thm2}, and Theorem \ref{conj11}}

\subsection{Proof of Theorem \ref{conj11} }
We first present the proof of Theorem \ref{conj11}  as it follows from a well-known trick.

\noindent
{\bf Proof of Theorem \ref{conj11}:}\label{pconj11} Note that $n=p(k-1)+q$. For the lower bound, let $G_3$ be the graph which consists of $p$ copies of $T(k-1,r-1)$ and one copy of $T(q,r-1)$. In addition, let $G_4=K_{\lfloor k/2 \rfloor-1}  \vee I_{n-\lfloor k/2 \rfloor+1}.$  Both $G_3$ and $G_4$ are $\{K_r,P_k\}$-free and the lower bound follows.

For the upper bound, assume that $G$ is a $\{K_r,P_k\}$-free graph with  $\lfloor k/2 \rfloor+1 \leq r <k  \leq n$. Let $u$ be a new vertex.  The graph resulted from connecting $u$  to each vertex of $G$ is denoted by $G'$. Notice that $G'$ contains $n'=n+1$ vertices and $e(G)+n$ edges. Additionally, $G'$ is $\{K_{r+1},\C_{\geq k+1}\}$-free. Let $r'=r+1$ and $k'=k+1$. Observe that $r' \geq \lfloor (k'-1)/2 \rfloor+2$, $r' <k' \leq n'$, and $n'-1=p(k'-2)+q$.

If $k$ is odd and $k=3$, then $G$ is a matching as $G$ is $P_3$-free which is trivial.  If $k$ is odd and $k \geq 5$, i.e., $k' \geq 6$ is even, then  Theorem \ref{thm1} yields that
\[
e(G') \leq p\ex(k'-1,K_{r'})+\ex(q+1,K_{r'}).
\]
As the assumption for $k$ and $r$, to get $T(k'-1,r'-1)$, one can add   a new vertex to $T(k-1,r-1)$ which is adjacent to all vertices in $T(k-1,r-1)$.  It follows that $\ex(k'-1,K_{r'})=\ex(k-1,K_{r})+k-1.$ Similarly, $\ex(q+1,K_{r'})=\ex(q,K_r)+q.$
Recall $n=p(k-1)+q$. Thus
\[
e(G)=e(G')-n  \leq p\ex(k'-1,K_{r'})+\ex(q+1,K_{r'})-n=p\ex(k-1,K_{r})+\ex(q,K_r).
\]
The case where $k$ is even can be proved similarly.  \hfill $\square$

\subsection{Proof of Theorem \ref{thm1}}
We sketch the proof of Theorem \ref{thm1} first. Let $G$ be a $\{K_r,\C_{\geq k}\}$-free graph. We apply the induction on the number of blocks in $G$. For even $k$, as we can show $\max \{g_r(n,(k-2)/2,k),g_r(n,2,k)\} \leq f(n,k,r)$ (see Lemma \ref{g2f} and Lemma \ref{gtf}), then the theorem follows from the subadditivity of $f(n,k,r)$ (see Lemma \ref{subadditivity}). For odd $k$, as
$g_r(n,(k-1)/2,k)>f(n,k,r)$ for large $n$, then the case where there is a block with $g_r(b,(k-1)/2,k)$ edges will be handled separately. The other case in which each block $B$ contains $g_r(b,2,k)$ edges also follows from  the subadditivity of $f(n,k,r)$.

We begin to prove a number of technical lemmas.  Readers can  skip proofs of lemmas first and go the proof of Theorem \ref{thm1} directly.

The next lemma gives a formula for the number of edges in a certain Tur\'an graph.
\begin{lemma} \label{edgesTuran}
	If  $r\ge\lfloor (k-1)/2\rfloor+2$ and $r-1 \le n\le k-1$, then $\ex(n,K_r)=\binom{n}{2}-(n-r+1)$.
\end{lemma}
{\bf Proof:} Notice that $\ex(n,K_r)$ is the number of edges in the Tur\'an graph $T(n,r-1)$. As the assumption $r \geq \lfloor (k-1)/2\rfloor+2 \geq k/2+1$ and $r \le n\le k-1$ , each partite of $T(n,r-1)$ contains at most two vertices and there are $n-r+1$ partites with exactly two vertices. Thus
\[
\ex(n,K_r)=\binom{n}{2}-(n-r+1),
\]
as required. \hfill$\square$

Recall the following definition.
For integers $n$ and $k$, if $n-1=p(k-2)+q$ with $q \leq k-3$, then
$$f(n,k,r)=p \ex(k-1,K_r)+\ex(q+1,K_r).$$
We next prove the subadditivity of the function $f(n,k,r)$.
\begin{lemma}\label{subadditivity}
	Let $n-1=p(k-2)+q$, $n_1-1=p_1(k-2)+q_1$ and $n_2-1=p_2(k-2)+q_2$, where $n_1+n_2-1=n$ and $0\le q,q_1,q_2\le k-3$. If $r \geq \lfloor (k-1)/2\rfloor+2$, then we have
	\[f(n_1,k,r)+f(n_2,k,r)\le f(n,k,r).\]
\end{lemma}
\noindent
{\bf Proof:} 
Let $G$ be the graph which consists of $p$ copies of $T(k-1,r-1)$ and one $T(q+1,r-1)$ all sharing a common vertex. Similarly, for $i=1,2$, let $G_i$ be the graph consists of $p_i$ copies of $T(k-1,r-1)$ and one $T(q_i+1,r-1)$ sharing a common vertex.
By the  definition, we have $f(n,k,r)=e(G)$ and $f(n_i,k,r)=e(G_i)$ for $i \in \{1,2\}$. As the assumption for $q_1$ and $q_2$,  it follows that either   $p_1+p_2=p$ or  $p_1+p_2=p+1$.

If $p=p_1+p_2$, then $q=q_1+q_2$ in this case. Removing edges contained in all $T(k-1,r-1)$ from $G$, $G_1$, and $G_2$, we can see that it is left to show $t(q+1,r-1) \geq t(q_1+1,r-1)+t(q_2+1,r-1)$ with the assumption $q=q_1+q_2.$ As the assumption
$r \geq \lfloor (k-1)/2\rfloor+2$ and $q_i \leq k-3$ for $i \in \{1,2\}$, there is at least one partite from $T(q_1+1,r-1)$ and $T(q_2+1,r-1)$ with one vertex. If one identify such a vertex from $T(q_1+1,r-1)$ with the one from $T(q_2+1,r-1)$, then the resulting graph is $K_r$-free  with $q_1+q_2+1$ vertices and $t(q_1+1,r-1)+t(q_2+1,r-1)$ edges.
 Now the inequality $t(q+1,r-1) \geq t(q_1+1,r-1)+t(q_2+1,r-1)$ follows from   the definition of the Tur\'an graph.

If $p=p_1+p_2+1$, then  $q_1+q_2=q+k-2$. It remains to establish $t(q_1+1,r-1)+t(q_2+1,r-1) \leq t(k-1,r-1)+t(q+1,r-1)$ under the condition $q_1+q_2=q+k-2$. Without loss of generality, we assume $q_1 \leq q_2$. For $2 \leq a \leq b \leq k-2$, we next show the fact that
$$t(a,r-1)+t(b,r-1) \leq t(a-1,r-1)+t(b+1,r-1).$$
As the assumption $r \geq \lfloor (t-1)/2\rfloor+2$, we observe that each partite in the Tur\'an graph $T(x,r-1)$ for $2 \leq x \leq k-1$ contains either two vertices or one vertex.
To see the fact, starting from the vertex disjoint Tur\'an graphs $T(a,r-1)$ and $T(b,r-1)$, we move one vertex $u$ from $T(a,r-1)$ to $T(b,r-1)$ such that the resulting graphs are $T(a-1,r-1)$ and $T(b+1,r-1)$.  By the observation above, through this process, we remove at most $a-1$ edges and create at least $b-1$ edges (the  vertex $u$ is contained in a partite with two vertices in $T(b+1,r-1)$). The fact follows easily as we assume $a \leq b$. We can apply this fact recursively staring with $a=q_1+1$, $b=q_2+1$, and  stop once one Tur\'an graph contains exactly $q+1$ vertices. \hfill $\square$

Recall  $g_r(n,a,k)=a(n-k+a)+\ex(k-a,K_r)$.  We first prove $g_r(n,2,k) \leq f(n,k,r)$.
 The following fact gives us the range in which $g_r(n_,2,k) > g_r(n,\lfloor (k-1)/2\rfloor,k)$.
\begin{lemma}\label{fact1}
Assume that $n \geq k \geq 7$ and $\lfloor (k-1)/2 \rfloor +2 \leq r <k $. If  $g_r(n_,2,k) > g_r(n,\lfloor (k-1)/2\rfloor,k)$, then $n \leq \tfrac{5k}{4}-1$.
\end{lemma}
{\bf Proof:} Suppose that $n>\tfrac{5k}{4}-1$. We next show
\[
K=g_r(n,\lfloor (k-1)/2 \rfloor,k)-g_r(n,2,k) \geq 0,
\]
 which is a contradiction to the assumption.
Recall that
$$g_r(n,a,k)=(n-k+a)a+\ex(k-a,K_r). $$
Moreover,  Lemma \ref{edgesTuran} yields that $\ex(k-2,K_r)=\binom{k-2}{2}-(k-r-1).$
Thus we have
\begin{align*}
g_r(n,2,k)&=2(n-k+2)+\ex(k-2,K_r)\\
&=2n+\frac{k^2}{2}-\frac{11}{2}k+r+8.
\end{align*}

If  $k$ is even, then we have  $r\ge k/2+1$. It implies that
\begin{align*}
	g_r\left(n,\frac{k-2}{2},k\right)&=\frac{k-2}{2}\left(n-\frac{k+2}{2}\right)+\ex\left(\frac{k+2}{2},K_r\right)\\
	&\ge \frac{k-2}{2}\left(n-\frac{k+2}{2}\right)+\binom{k/2+1}{2}-1\\
	&=\frac{k-2}{2}n-\frac{k^2}{8}+\frac{k}{4}.
\end{align*}
Recall assumptions $n>\tfrac{5k}{4}-1$, $k \geq 8$ ($k$ is even and $k \geq 7$), and $r<k$. Thus
\begin{align*}
	K&\ge \frac{k-2}{2}n-\frac{k^2}{8}+\frac{k}{4}-\left(2n+\frac{k^2}{2}-\frac{11}{2}k+r+8\right)\\
	&=\left(\frac{k}{2}-3\right)n-\frac{5}{8}k^2+\frac{23}{4}k-r-8\\
	&>\left(\frac{k}{2}-3\right)\left(\frac{5}{4}k-1\right)-\frac{5}{8}k^2+\frac{23}{4}k-r-8\\
	&=\frac{3k}{2}-r-5\\
	&\ge\frac{3k}{2}-(k-1)-5\\
    &=k/2-4\ge 0,\\
\end{align*}
as desired.

If $k$ is odd, then we have $r\ge\frac{k+3}{2}.$ It follows that
\begin{align*}
	g_r\left(n,\frac{k-1}{2},k \right)&=\frac{k-1}{2}\left(n-\frac{k+1}{2}\right)+\ex\left(\frac{k+1}{2},K_r\right)\\
	&= \frac{k-1}{2}\left(n-\frac{k+1}{2}\right)+\binom{(k+1)/2}{2}\\
	&=\frac{k-1}{2}n-\frac{k^2}{8}+\frac{1}{8}.
\end{align*}
Therefore,
\begin{align*}
	K&=\frac{k-1}{2}n-\frac{k^2}{8}+\frac{1}{8}-\left(2n+\frac{k^2}{2}-\frac{11}{2}k+r+8\right)\\
	  &=\frac{k-5}{2}n-\frac{5}{8}k^2+\frac{11}{2}k-r-\frac{63}{8}\\
	&>\frac{k-5}{2}\left(\frac{5k}{4}-1\right)-\frac{5}{8}k^2+\frac{11}{2}k-r-\frac{63}{8}\\
	&=\frac{15k}{8}-\frac{43}{8}-r\\
	&\ge\frac{15k}{8}-\frac{43}{8}-(k-1)\\
	&=\frac{7}{8}(k-5)>0,
\end{align*}
as desired. The proof is complete. \hfill $\square$

We next show $g_r(n,2,k) \leq  f(n,k,r)$.
\begin{lemma}\label{g2f}
	If $6 \leq k\le n\le \tfrac{5k}{4}-1$,  $r \geq \lfloor (k-1)/2 \rfloor +2$, and   $n-1=k-2+q$ with $q \leq k-3$,  then
	\[
	g_r(n,2,k)\le f(n,k,r).
	\]
\end{lemma}
\noindent
{\bf Proof:}  As $ k \le n\le  \tfrac{5k}{4}-1$, it follows that  $q \le \tfrac{k}{4}.$  Since $r\ge \lfloor (k-1)/2\rfloor+2 \ge k/2+1,$  we have $r-1\ge q+1$ and $\ex(q+1,K_r)=\binom{q+1}{2}.$  By Lemma \ref{edgesTuran}, we have
\[
\ex(k-1,K_r)=\binom{k-1}{2}-(k-r).
\]
Similarly,  $\ex(k-2,K_r)=\binom{k-2}{2}-(k-r-1)$.  Then
\begin{align*}
    M&=f(n,k,r)- g_r(n,2,k)\\
	&=\binom{k-1}{2}-(k-r)+\binom{q+1}{2}-2(n-k+2)-\binom{k-2}{2}+(k-r-1)\\
	&=k-3+\binom{q+1}{2}-2(n-k+2).
\end{align*}
If $q \geq 2$, then $M \geq k-2(n-k+2) \geq \tfrac{k}{2}-2 \geq 0$ as  $n \leq \tfrac{5k}{4}-1$ and $k \geq 6$.  As $n \geq k$, it is left to consider the case where $q=1$ and $n=k$. Note that $M \geq 0$ in this case since the assumption $k \geq 6$.
\hfill $\square$

For even $k$, the next lemma shows $g_r(n,(k-2)/2,k)\le f(n,k,r)$.
\begin{lemma}\label{gtf}
	Let $n-1=p(k-2)+q$ with $q \leq k-3$.  If $k$ is even,  $n\ge k\ge 6$, and $k/2+1 \leq r \leq k-1$, then $g_r(n,(k-2)/2,k)\le f(n,k,r)$.
\end{lemma}
\noindent
{\bf Proof:}
Note that $(k-2)/2=\lfloor (k-1)/2\rfloor$ as $k$ is even.
Additionally,
\[
g_r(n,(k-2)/2,k)=\frac{k-2}{2}(n-(k+2)/2)+\ex((k+2)/2,K_r).
\]
Notice that $\ex(k-1,K_r)$ is the number of edges in the Tur\'an graph $T(k-1,r-1)$. As the assumption $k/2+1 \leq r \leq k-1$,  by Lemma \ref{edgesTuran}, we have
\[
\ex(k-1,K_r)=\binom{k-1}{2}-(k-r).
\]
Similarly, $\ex(q+1,K_r)=\binom{q+1}{2}-(q-r+2)$ for $q+1 \geq r$ and $\ex(q+1,K_r)=\binom{q+1}{2}$ for $q+1 \leq r-1$.
Recall
\[
f(n,k,r)=p\ex(k-1,K_r)+\ex(q+1,K_r).
\]
\noindent
{\bf Case 1:} $r \geq k/2+2$. Apparently,
\[
\ex((k+2)/2,K_r)=\binom{(k+2)/2}{2}  \textrm{ and }  g_r(n,(k-2)/2,k)=\frac{(k-2)n}{2}-\frac{k^2}{8}+\frac{k}{4}+1.
\]
If $q+1 \geq r$, then
\begin{align*}
	f(n,k,r)&=p\binom{k-1}{2}-p(k-r)+\binom{q+1}{2}-(q-r+2)\\
	&=\frac{n-1-q}{k-2}\binom{k-1}{2}+\binom{q+1}{2}-\frac{n-1-q}{k-2}(k-r)-(q-r+2)\\
	&=\frac{(k-1)n}{2}-\frac{(q+1)(k-1)}{2}+\frac{q(q+1)}{2}-\frac{n-1-q}{k-2}(k-r)-(q-r+2)\\
	& \geq \frac{(k-1)n}{2}-\frac{(q+1)(k-1)}{2}+\frac{q(q+1)}{2}-\frac{n-1-q}{k-2}(\frac{k}{2}-2)-(q-r+2) \hfil \textrm{ (as } r \geq k/2+2) \\
	&=\frac{(k-2)n}{2}+\frac{(q+1)^2}{2}-\frac{(q+1)(k+1)}{2}+\frac{n-1-q}{k-2}+r-1 \\
	& \geq  \frac{(k-2)n}{2}-\frac{k^2}{8}-\frac{k}{4}+1+\frac{n-1-q}{k-2}+r-1
	\hfil ( \textrm{ as the minimum is achievd for } q+1=r/2+2)\\
	& \geq \frac{(k-2)n}{2}-\frac{k^2}{8}+\frac{k}{4}+1\\
	&=g_r(n,(k-2)/2,k).
\end{align*}
For $q+1 \leq r-1$, we get that
\begin{align*}
f(n,k,r)&=p\binom{k-1}{2}-p(k-r)+\binom{q+1}{2}\\
	&=\frac{n-1-q}{k-2}\binom{k-1}{2}+\binom{q+1}{2}-\frac{n-1-q}{k-2}(k-r)\\
	&=\frac{(k-1)n}{2}-\frac{(q+1)(k-1)}{2}+\frac{q(q+1)}{2}-\frac{n-1-q}{k-2}(k-r)\\
	&\geq \frac{(k-1)n}{2}-\frac{(q+1)(k-1)}{2}+\frac{q(q+1)}{2}-\frac{n-1-q}{k-2}(\frac{k}{2}-2) \hfil \textrm{ (as } r\geq k/2+2) \\
	&=\frac{(k-2)n}{2}+\frac{(q+1)^2}{2}-\frac{(q+1)(k-1)}{2}+\frac{n-1-q}{k-2} \\
	& \geq  \frac{(k-2)n}{2}-\frac{k^2}{8}+\frac{k}{4}+\frac{n-1-q}{k-2} \hfil (\textrm{ as the minimum is achievd for } q+1 \in \{k/2,k/2-1\})\\
	& \geq \frac{(k-2)n}{2}-\frac{k^2}{8}+\frac{k}{4}+1 \hfil (\textrm{ as } p=\frac{n-1-q}{k-2} \geq 1)\\
	&=g_r(n,(k-2)/2,k).
\end{align*}
{\bf Case 2:} $r=k/2+1$. Observe that
\[
\ex((k+2)/2,K_r)=\binom{(k+2)/2}{2}-1  \textrm{ and }  g_r(n,k/2-1,k)=\frac{(k-2)n}{2}-\frac{k^2}{8}+\frac{k}{4}.
\]
Repeating the argument in Case 1, we can show $$f(n,k,r) \geq g_r(n,k/2-1,k)$$ similarly.
\hfill $\square$

As we did before, let $t=\lfloor (k-1)/2\rfloor$. Recall
$$g_r(n,a,k)=(n-k+a)a+\ex(k-a,K_r)$$
 and
 $$f(n,k,r)=p \ex(k-1,K_r)+\ex(q+1,K_r).$$

\noindent
{\bf Proof of Theorem \ref{thm1}: even $k$:}  To see the lower bound, observe that $F(n,k,r)$  is $\{K_r,\C_{\geq k}\}$-free and $e(F(n,k,r))=f(n,k,r)$. The lower bound follows.


For the upper bound, let $G$ be an $n$-vertex $\{K_r,\C_{\geq k}\}$-free graph with maximum number of edges. We claim that $G$ is connected. Otherwise, if $G$ is not connected, then the graph resulted from adding an edge to connect two connected components is still $\{K_r,\C_{\geq k}\}$-free as both $K_r$ and $\C_{\geq k}$ are 2-connected, but has more edges than $G$. This is a contradiction to the definition of $G$. Thus $G$ is connected. We apply the induction on the number of blocks contained in $G$ to show the upper bound. If $G$ contains one block ($G$ is 2-connected), then $e(G) \leq \ex(n,K_r)=f(n,k,r)$ for $n \leq k-1$ and
$$
e(G) \leq \max\{g_r(n,2,k),g_r(n,t,k)\} \leq f(n,k,r)
$$
for $n \geq k$. The inequality above follows from the combination of Theorem \ref{thm11},  Lemma \ref{g2f}, and Lemma \ref{gtf}.

Assume that $G$ contains at least two blocks. Let $B$ be an end block of $G$ and $b$ be the cut vertex. Let $G'=(G-B)\cup\{b\}$ and $n'=|V(G')|$. Then $G'$ contains fewer blocks than $G$.  If $n' <k$, then $e(G') \leq \ex(n,K_r)=f(n',k,r)$. Otherwise, by induction hypothesis, $e(G') \leq f(n',k,r)$. Additionally,  $e(B) \leq f(n-n'+1,k,r)$ by using the proof for the base case.
Therefore, by Lemma \ref{subadditivity},
\begin{align*}
	e(G)=e(G')+e(B)
	&\le f(n',k,r)+f(n-n'+1,k,r)\\
	&\leq f(n,k,r).
\end{align*}
The proof of the even case is complete.
\hfill $\square$

\vspace{0.1cm}
\noindent
{\bf Proof of Theorem \ref{thm1}: odd $k$:} For the lower bound,  we already showed that $F(n,k,r)$ is $\{K_r,\C_{\geq k}\}$-free. Note that  the clique number of $G_1$ is $t+1$ which is less that $r$ by the assumption of $r$. In addition, a longest cycle contains at most $k-1$ edges. Thus $G_1$ is also $\{K_r,\C_{\geq k}\}$-free.  The lower bound is proved.

For the upper bound, we assume $k \geq 7$ for a moment. As we did above, let $G$ be an $n$-vertex $\{K_r,\C_{\geq k}\}$-free graph with maximum number of edges. Then $G$ is connected. Since both $K_r$ and $\C_{\geq k}$ are 2-connected, we can further assume that  the subgraph induced by each block contains the maximum number of edges. In other words, if $B$ is a block and $b=|V(B)|$, then either $e(G[B])=\ex(b,K_r)$ for $2\le b\le k-1$ or $e(G[B])=\max\{g_r(b,2,k),g_r(b,t ,k)\}$ for $b\ge k$.
Apparently,   each block $B$ is $\C_{\geq k}$-free.  Note that $(k-1)/2=\lfloor (k-1)/2\rfloor=t$ as $k$ is odd. Theorem \ref{EGcycle} yields that  $e(B) \leq t(|V(B)|-1)$ no matter the size of $B$.

Assume $B_1,\ldots,B_j$ are blocks of $G$. For each $1 \leq i \leq j$, let $b_i$ be the number of vertices in $B_i$.
Consider a rooted block tree with $B_1$ as the root. Observe that each block $B_i$ brings exactly $b_i-1$ new vertices and then $\sum_{i=2}^j (b_i-1)=n-b_1$. For the case where $G$ contains a block, say $B_1$, such that $e(B_1)=g_r(n,t,k)$, then
\begin{align*}
	e(G)=&\sum_{i=1}^{j}e(B_i)=e(B_1)+\sum_{i=2}^{j} e(B_i) \\
	&\leq e(B_1)+ \sum_{i=2}^{j} t(b_i-1) \\
	&=g_r(b_1,t,k) + \sum_{i=2}^{j}t(b_i-1)\\
	&=g_r(b_1,t,k)+ t(n-b_1)\\
	&=\binom{t}{2}+t(b_1-1)-t(t-1)+t(n-b_1)\\
	&=\binom{t}{2}+t(n-1)-t(t-1)\\
	&=g_r(n,t,k).
\end{align*}
It remains to consider the case where each block $B_i$ in $G$ satisfies $e(B_i) \neq g_r(n,t,k)$. We apply the induction on the number of blocks in $G$ to show $e(G) \leq f(n,k,r)$. For the case where $G$ contains one block (i.e., $G$ is 2-connected), if $n \geq k$, then $e(G) = g_r(n,2,k) \leq f(n,k,r)$ by the assumption of $G$ and Lemma \ref{g2f}.
If $n <k$, then $e(G) \leq \ex(n,K_r)=f(n,k,r)$.
For the case where $G$ contains at least two blocks,
  let $B$ be an end block of $G$ and $b$ be the cut vertex. Set $G'=(G-B)\cup\{b\}$. Then $G'$ contains fewer blocks than $G$. Let $n'=|V(G')|$. Then $e(G') \leq f(n',k,r)$ for $n' \geq k$ by the induction hypothesis and $e(G') \leq \ex(n,K_r)=f(n',k,r)$ for $n' <k$.
   Moreover,  repeating  the proof for base case, we can show $e(G[B])\le f(b,k,r)=f(n-n'+1,k,r)$. 
By induction hypothesis and Lemma \ref{subadditivity}, it follows that
\begin{align*}
	e(G)=e(G')+e(B)
	&\le f(n',k,r)+f(n-n'+1,k,r)\\
	&\leq f(n,k,r).
\end{align*}
The upper bound is proved for $k \geq 7$. It remains to consider the case where $k=5$. In this case, note that $r=4$ and $e(B_i)=g_4(b_i,2,5)$ whenever $b_i \geq 5$. Moreover, $e(B_i) \leq 2(b_i-1)$ still holds.
 We can repeat  the argument for  $k \geq 7$ to show the upper bound.
 The proof for Theorem \ref{thm1} is complete.
 \hfill $\square$

\subsection{Proof of Theorem \ref{thm2}}
To prove Theorem \ref{thm2}, let $G$ be a $\{K_r,\C_{\geq k}\}$-free graph. We also apply the induction on the number blocks contained in $G$. Since Theorem \ref{thm2} only works for large $n$, we need to establish an upper bound for the number of edges in small blocks. The proof will be split into two cases depending on whether
there is a large block.

We again fix $t=\lfloor (k-1)/2\rfloor$.
To prove the upper bound in Theorem \ref{thm2}, we need two more specific lemmas.
\begin{lemma} \label{onemore}
	If $2\le  n<k,k\ge 5$ and $3\le r\le t+1$, then $\ex(n,K_r)\le (t-1/k)(n-1).$
\end{lemma}
{\bf Proof:} 	For the case that  $n\le k-3$, we have
\[
\ex(n,K_r)\le \binom{n}{2}=\frac{n(n-1)}{2}\le \frac{(k-3)(n-1)}{2}\le(t-\frac{1}{k})(n-1)
\]
as $t=\lfloor(k-1)/2\rfloor.$
It remains  to show that the inequality holds for $n=\{k-2,k-1\}.$
If $n=k-2$ and $k$ is odd, then we also have $\ex(n,K_r)\le(t-\frac{1}{k})(n-1)$ by the argument above. If $n=k-2$ and $k$ is even, then $r-1\le t= (k-2)/2$ and each partite of $T(k-2,r-1)$ contains at least two vertices. As $r-1 \geq 2$ and $1-3/k <1$, we have
\[
 \ex(k-2,K_r)\le\frac{(k-2)(k-3)}{2}-(r-1)\le\frac{k-2}{2}(k-3)-(1-3/k)=\left(t-\frac{1}{k}\right)(k-3).
\]
For the case that $n=k-1$, as $r-1\le t$,  we have
\[\ex(k-1,K_r)\le \left(1-\frac{1}{r-1}\right)\frac{(k-1)^2}{2}\le \left(1-\frac{1}{t} \right)\frac{(k-1)^2}{2}\le \left(t-\frac{1}{k} \right)(k-2).\]
The last inequality can be verified directly by noting that $t=\lfloor (k-1)/2\rfloor$ and $k\ge 5.$ \hfill$\square$
\begin{lemma} \label{twomore}
	If $5 \le k\le n$ and $3\le r\le t+1$, then
	\[	
	\ex(k-1,K_r)+\left(t-\frac{1}{2} \right)(n-k+1)\le \left(t-\frac{1}{k} \right)(n-1).
	\]		
\end{lemma}
{\bf Proof:} Lemma \ref{onemore} gives that $\ex(k-1,K_r) \leq (t-\tfrac{1}{k})(k-2)$. Apparently, $(t-\tfrac{1}{2})(n-k+1)\le (t-\tfrac{1}{k})(n-k+1)$ as $k \geq 5$. The lemma follows.  \hfill$\square$

We are now ready to prove Theorem \ref{thm2}.\\
\noindent
{\bf Proof of Theorem \ref{thm2}:} For the lower bound, as the clique number of   $G_2$  is $r-1$ and a longest cycle contains at most $2t \leq k-1$ edges. Thus $G_2$ is $\{K_r,\C_{\geq k}\}$-free and  the lower bound follows.

To make the notation simple, we define $$g'_r(n,t,k)=(n-t)t+\ex(t,K_{r-1}).$$
Note that
$$g'_r(n,t,k)=t(n-1)-t(t-1)+\ex(t,K_{r-1}) \leq t(n-1).$$

 We next prove the upper bound.  Let $G$ be an $n$-vertex $\{K_r,\C_{\geq k}\}$-free graph with maximum number of edges, where $n \geq \tfrac{k^3}{4}$. As we did in the proof of Theorem \ref{thm1}, we can assume that $G$ is connected. Let $B_1,\ldots,B_j$ be blocks of $G$. For $1 \leq i \leq j$, let $b_i$ be the size of $B_i$.

If $b_i \geq \tfrac{k^2}{2}$, then Theorem \ref{thm21} gives us  $e(B_i) \le t(b_i-t)+\ex(t,K_{r-1})=g'_r(b_i,t,k).$
If $k \leq b_i < \tfrac{k^2}{2}$, then by Lemma \ref{newstability}, either $e(B_i) \leq t(b_i-t)+\ex(t,K_{r-1})=g'_r(b_i,t,k)$
 or
$$
e(B_i) \leq \max\{(t-1/2) (b_i-1),\ex(k-1,K_r)+(t-\tfrac{1}{2})(b_i-k+1)\} \leq (t-1/k)(b_i-1).
$$
The last inequality follows from Lemma \ref{twomore}. If $b_i<k$, then we also have
\[
e(B_i)\le\ex(b_i,K_r)\le(t-1/k)(b_i-1).
\]
by Lemma \ref{onemore}. Therefore, for any $b_i <\tfrac{k^2}{2}$, we have either $e(B_i) \leq g'_r(b_i,t,k)$ or $e(B_i)  \leq (t-1/k)(b_i-1)$. Note that $e(B_i) \leq  t (b_i-1)$ holds for any $1 \leq i \leq j$ as $B_i$ is $\C_{\geq k}$-free and Theorem \ref{EGcycle}.

For the case where there is a $B_i$, say $B_1$,  such that $b_1 \geq \tfrac{k^2}{2}$, then $e(B_1) \leq g'_r(b_1,t,k)$.  Note that $\sum_{i=2}^j (b_i-1)=n-b_1$ by considering a rooted block tree with $B_1$ as the root. We can upper bound $e(G)$ as follows:
\begin{align*}
	e(G)&=\sum_{i=1}^{j}e(B_i)\\
         &=e(B_1)+\sum_{i=2}^{j} e(B_i) \\
	&\leq e(B_1)+ \sum_{i=2}^{j} t(b_i-1) \\
	&=g_r'(b_1,t,k) + \sum_{i=2}^{j} t(b_i-1)\\
	&=g_r'(b_1,t,k)+ t(n-b_1)\\
	&=t(b_1-1)-t(t-1)+\ex(t,K_{r-1})+t(n-b_1)\\
	&=t(n-1)-t(t-1)+\ex(t,K_{r-1})\\
    &=t(n-t)+\ex(t,K_{r-1}).
\end{align*}
It remains to consider the case where $b_i<\tfrac{k^2}{2}$ for each $1 \leq i \leq j$. If there is a $B_i$ such that
$e(B_i) \leq g'_r(b_i,t,k)$, then we can repeat the argument above to show the desired upper bound.  We are left to show the upper bound for the case where $b_i<\tfrac{k^2}{2}$ and $e(B_i) \leq (t-1/k)(b_i-1)$ for each $1 \leq i \leq j$. For this case, we have
\begin{align*}
	e(G)&=\sum_{i=1}^{j}e(B_i)\\
	& \leq \sum_{i=1}^j (t-1/k)(b_i-1)\\
	&=(t-1/k)(n-1)\\
	& \leq t(n-1)-t(t-1)+\ex(t,K_{r-1})\\
   &=t(n-t)+\ex(t,K_{r-1}),
\end{align*}
the second-to-last inequality follows from the assumption $n \geq \tfrac{k^3}{4}$.  The proof is complete. \hfill $\square$
\section{Concluding remarks}
In this paper, we proved results on the Tur\'an number of $\h=\{K_r,\C_{\geq k}\}$.  For $r \geq \lfloor (k-1)/2 \rfloor+2$, we obtained the value of $\ex(n,\h)$ for all $n$. However, for $3 \leq r \leq \lfloor (k-1)/2 \rfloor+1$, we were only able to show the value of $\ex(n,\h)$ for large $n$. A natural question is to determine $\ex(n,\h)$ for small $n$ in this case.
Let us recall
lower bound constructions in this paper (Theorem \ref{thm1} and Theorem \ref{thm2}).  For  the graph $F(n,k,r)$,
we  start with a $\C_{\geq k}$-free graph with the maximum number of edges and turn it  $K_r$-free. To construct $G_2$,
we modify the Tur\'an graph $T(n,r-1)$ so that the sum of size of smallest $(r-2)$-partite is at most $\lfloor (k-1)/2 \rfloor$ (to ensure it is $\C_{\geq k}$-free). The answer to the question above requires new ideas to construct $\h$-free graphs.

 Note that upper bounds for $\ex(n,P_k)$ and $\ex(n,\C_{\geq k})$ in Theorem \ref{EGpath} and Theorem \ref{EGcycle} are only tight for particular $n$. However, these upper bounds are widely applicable because they are very simple. Therefore, it is desired to show such kind upper bound for $\ex(n, \{K_r,\C_{\geq k}\})$.  Unfortunately, this kind of upper bound does not exist in general as Theorem \ref{thm1} involves two extremal graphs for odd $k$.
Therefore, in comparison to the Tur\'an function $\ex(n,\C_{\geq k})$, the Tur\'an function $\ex(n,\{K_r,C_{\ge k}\})$ indeed behaves differently.

As we mentioned before,  the result on the Tur\'an number of $P_k$ is a simple corollary of the one for $\C_{\geq k}$. One may ask whether it is still the case for $K_r$-free graphs.  The answer is positive for $r \geq \lfloor (k-1)/2 \rfloor+2$ as Theorem \ref{conj11}. For $r \leq \lfloor (k-1)/2 \rfloor+1$, we are not able to prove such kind result. One possible reason is that $\ex(n, \{K_r,\C_{\geq k}\})$ is not known for small $n$. If one can determine $\ex(n, \{K_r,\C_{\geq k}\})$ for all $n$, then it is very likely to provide a positive answer to the question mentioned above.

  It is obvious that $\ex_{2\textrm{-conn}}(n,\h) \leq \ex(n,\h)$ for any $\h$ by  definitions. Let  $\h=\{K_r,\C_{\geq k}\}$ and $r=k/2+1$ with $k$ being even.  If $n-1=p(k-2)+q$ with $p \geq 1$ and $q \in \{k/2-1,k/2-2\}$, then   $\max\{g_r(n,2,k),g_r(n,k/2-1,k)\}=g_r(n,k/2-1,k)= f(n,k,r)$ by the direct computation.
  Therefore, $\ex_{2\textrm{-conn}}(n,\h)= \ex(n,\h)$ by Theorem \ref{thm1} and Theorem \ref{thm11}, which indicates that
 $r=\lfloor (k-1)/2 \rfloor+2$ is indeed a critical case.

\end{document}